\documentclass[12pt]{article}   
\usepackage{latexsym}        
\usepackage{amssymb}         
\usepackage{amscd}           

\usepackage{amssymb,amsthm,amsmath}

\input epsf

\newtheorem{defn0}{Definition}[section]

\newtheorem{thm0}[defn0]{Theorem}
\newtheorem{lemma0}[defn0]{Lemma}
\newtheorem{corollary0}[defn0]{Corollary}

\newtheorem{remark0}[defn0]{Remark}

\newtheorem{nota0}[defn0]{Notation}

\newcommand{\defref}[1]{Definition~\ref{#1}}

\newcommand{\thmref}[1]{Theorem~\ref{#1}}

\makeindex

\font\Bbb=msbm10

\newcommand{\HHH}{\mbox{\Bbb H}}

\newcommand{\PP}{\mbox{\Bbb P}}

\newcommand{\CC}{\mbox{\Bbb C}}

\newcommand{\OO}{\mbox{${\cal O}$}}

 \def\rig#1{\smash{ \mathop{\longrightarrow}
    \limits^{#1}}}

\begin{document}

\begin{center}
 {\Large \bf NONDEGENERATE MULTIDIMENSIONAL MATRICES AND INSTANTON BUNDLES
}
\end{center}

\begin{center}
LAURA COSTA \footnote{Partially supported by DGICYT
PB97-0893.} \hspace{8mm} GIORGIO OTTAVIANI
\footnote{Partially supported by Italian MURST.} \end{center}

\centerline{November 27, 2000}

\vspace{3mm}


{\bf Abstract: }In this note we prove that the moduli space of rank $2n$
symplectic instanton bundles on
${\PP^{2n+1}}$,  defined from the well known monad condition,
is affine. This result was not known even in the case $n=1$,
where by
\cite{ADHM} the real instanton bundles correspond to self dual Yang Mills
$Sp(1)$-connections over
 the $4$-dimensional sphere.
The result is proved as a consequence of the existence of an invariant of the
multidimensional matrices representing
the instanton bundles.

\section{Introduction} \indent

 A symplectic instanton bundle on $\PP^{2n+1}_{\CC }$ is a bundle of rank $2n$
defined as the
cohomology bundle of a well known monad (see \defref{def:instanton}).

In \cite{ADHM} it was shown that  instanton bundles on $\PP^3$ satisfying a
reality
condition correspond to  self dual  Yang Mills $Sp(1)$-connections over the
$4$-dimensional sphere $S^4=\PP^1_{\HHH}$. This correspondence was
generalized by
Salamon (\cite{Sal84}) who showed that instanton bundles on $\PP^{2n+1}$
which are trivial on the fiber of the twistor map $\PP^{2n+1}\to\PP^n_{{\HHH
}}$
correspond to $Sp(n)$-connections which minimize a certain Yang Mills
functional over
$\PP^n_{\HHH}$. We denote by $MI_{\PP^{2n+1}}(k)$  the moduli space of
 symplectic instanton bundles on $\PP^{2n+1}$ with $c_2=k$ (see
\defref{gitquot}).
Few things are known on $MI_{\PP^{2n+1}}(k)$  in general. It is known
(\cite{ADHM})
that it has a component of dimension $8k-3$ for $n=1$,
 it is smooth for $n=1$ and $k\le 5$ (\cite{KO99})
and it is expected to be singular and reducible for $n\ge 2$
and $k\ge 4$ (see \cite{AO00}).

Our main result, proved in section 2, is the following

\begin{thm0}
\label{affine0}
 $MI_{\PP^{2n+1}}(k)$ is affine.
\end{thm0}

In particular  $MI_{\PP^{2n+1}}(k)$  does not contain any complete subvariety
of
positive dimension. For other results of this type see \cite{HH}.
The technique we use to
prove the main theorem is to exhibit $MI_{\PP^{2n+1}}(k)$
as the GIT-quotient of an affine variety ${\cal Q}^0$ and then use standard
results
about invariant theory. The fact that ${\cal Q}^0$ is affine is a consequence
of
the existence of a invariant of multidimensional matrices representing the
instanton bundles, which generalizes the hyperdeterminant (see \cite{GKZ}
and \cite{AO99}).

\vspace{3mm}

The first named author would like to thank the Dipartimento di Matematica U.
 Dini for their hospitality and support at the time of the preparation of this
 paper.

\section{Notations and preliminaries} \indent

 We will start fixing some notation and recalling some facts about
$k$-instanton
 bundles on $\PP^{2n+1}=\PP(V)$, where $V$ is a complex vector space of
dimension $2n+2$ (See for instance \cite{OS86} and \cite{AO94}).

\begin{nota0}
\rm $\OO(d)=\OO_{\PP^{2n+1}}(d)$ denotes the invertible sheaf of degree
$d$ on $\PP^{2n+1}$ and for any coherent sheaf $E$ on $\PP^{2n+1}$ we
denote $E(d)=E \otimes \OO_{\PP^{2n+1}}(d)$ .
\end{nota0}

\begin{defn0}
\label{def:instanton} A symplectic instanton bundle $E$ over
$\PP^{2n+1}=\PP(V)$
is a bundle of rank $2n$ which appears as a cohomology bundle of a
monad
\refstepcounter{equation}
 \[ (\theequation) \label{definstanton} \hspace{8mm}
 I^*\otimes\OO(-1)\rig{A} W\otimes\OO\rig{A^t} I\otimes\OO(1) \]
where $(W,J)$ is a symplectic complex vector space of
dimension $2n+2k$ and $I$ is a complex vector space of
dimension $k$.
\end{defn0}

We do not assume in the definition that $E$ is stable, so we have to recall
some results.

The monad condition means that $A$ is injective (as a bundle morphism), $A^t$
is surjective and $ \mathrm{im } A\subset\ker A^t$ so that $E\simeq
\ker A^t/\mathrm{im } A$.
The fact that the map
\break\noindent $W\otimes\OO\rig{A^t} I\otimes\OO(1)$ is surjective, is
equivalent
to the fact that
the matrix $A\in Hom(V^*\otimes I^*,W)$ representing $E$
is nondegenerate according to \cite{GKZ}
(see \defref{def:dege} for the precise definition).

$Hom(V^*\otimes I^*,W)$ contains the subvariety ${\cal Q}$
 given by matrices $A$
for which the sequence (\ref{definstanton}) is a complex,
 that is, such that $A^tJA=0$.
$GL(I)\times Sp(W)$ acts on ${\cal Q}$ by
$(g,s)\cdot A=sAg$.

\begin{defn0}
\label{def:dege}
A matrix $A\in Hom(V^*\otimes I^*,W)$ is called degenerate
if the multilinear system $A(v\otimes i)=0$ has a solution such that
 $0\neq v\in V^*$ and $0\neq i\in I^*.$
\end{defn0}

By \cite{GKZ}, Theorem 14.3.1, this is equivalent to the standard
definition of degeneracy given in chapter 14.1 of \cite{GKZ}.
It is easy to check that nondegenerate matrices fill an irreducible
 subvariety $N$ of  $Hom(V^* \otimes I^*,W)$
of codimension $k$ (see \cite{WZ}). Hence only in the case $k=1$ it
 is well defined a hyperdeterminant according to
\cite{GKZ}. In the next section we will define a
$SL(I)\times Sp(W)$-invariant on $Hom(V^*\otimes
I^*,W)$ called $D$ which generalizes the hyperdeterminant and it is suitable
for our purposes.

 It was shown in
\cite{AO94} that all instanton bundles are simple, so that they
carry a unique symplectic form. Moreover for $n=1, 2$
it was proved in \cite{AO94} that all instanton bundles are stable,
and it is expected that the same result is true for $n\ge 3$.
In \cite{BH} it was essentially proved that there is a natural 1:1
correspondence between

\noindent {i)} isomorphism classes of symplectic instanton bundles, and

\noindent{ii)} orbits of $GL(I)\times Sp(W)$ on the open subvariety ${{\cal
Q}^0}$ of ${\cal Q}$ given by nondegenerate matrices.

In fact, using the quoted results of
\cite{AO94}, one can see that $\cite{BH}$; Section $4$ and Theorem at page 19,
adapt literally to our situation.
Moreover, \cite{BH}; Remark 2 at page 19  shows that, if $G$
denotes the quotient of $GL(I)\times Sp(W)$ by $\pm (id, id)$, then
 $G$ acts freely on ${{\cal Q}^0}$. In particular all points of
 ${{\cal Q}^0}$ are stable for the action of $GL(I)\times Sp(W)$.

\begin{defn0}
\label{gitquot}
The GIT-quotient ${{\cal Q}^0}/{GL(I)\times Sp(W)}$
is denoted by $MI_{\PP^{2n+1}}(k)$ and it is called
the moduli space of  symplectic  $k$-instanton bundles on $\PP^{2n+1}$.
It is a geometric quotient.
\end{defn0}

The above discussion shows that $MI_{\PP^{2n+1}}(k)$ coincides for $n=1,2$
with the open subset ${\cal MI}_{\PP^{2n+1}}(k)$ of the Maruyama scheme of
symplectic stable bundles on $\PP^{2n+1}$ of rank $2n$ and Chern polynomial
$\frac{1}{(1-t^2)^k}$ which are instanton bundles (this is an open condition
because by Beilinson theorem is equivalent to certain vanishing in cohomology,
see \cite{OS86} ).
In particular our notation for $MI_{\PP^{3}}(k)$ is consistent
with the usual one.  For $n\ge 3$ it is expected that the same result is true,
 but at present we can only say that ${\cal MI}_{\PP^{2n+1}}(k)$ is
 an open subset of $MI_{\PP^{2n+1}}(k)$.

\vspace{4mm}
\section{The invariant $D$ and the proof of the main result} \indent

First of all we remark that the vector spaces
$W\otimes S^nI$ and $V\otimes S^{n+1}I$ have the same dimension
$(2n+2k){\binom{k+n-1}{n}}=(2n+2){\binom{k+n}{n+1}}$.
We can construct from
\[W\rig{A^t} V\otimes I\]
the following morphisms
\refstepcounter{equation}
\[ \begin{array}{c} A^t\otimes id_{S^nI}:  W\otimes S^nI\rig{} V\otimes
I\otimes
S^nI ,\\
  id_V\otimes\pi : V\otimes I\otimes S^nI\rig{} V\otimes S^{n+1}I,
\end{array}\]
where $\pi$ is the natural projection, and
we consider the composition
 \[ (\theequation) \label{magic} \hspace{8mm}
   \Delta_A=(id_V\otimes\pi)\cdot(A^t\otimes id_{S^nI})
\colon W\otimes S^nI\rig{}V\otimes S^{n+1}I.  \]
\begin{defn0}
\label{def:magic} Let $A\in Hom(V^*\otimes I^*,W)$. We define $D(A)$ to be the
usual determinant
of the morphism $\Delta_A$ in (\ref{magic}) induced by $A$.
\end{defn0}

Notice that
\[ D\colon Hom(V^*\otimes I^*,W)\to (\det I)^{\alpha}\otimes (\det V)^{\beta}\]
where $\alpha=2\binom{k+n}{n}$ and $\beta=\binom{k+n}{n+1}$
is a $GL(V)\times GL(I)\times Sp(W)$-equivariant map
and $D(A)=0$ defines a homogeneous hypersurface of degree
$(2n+2k){\binom{k+n-1}
{n}}=(2n+2){\binom{k+n}{n+1}}$.
After a basis has been fixed in each of the vector spaces $V$, $I$ and $W$,
the map $D$ can be seen as a  $SL(V)\times SL(I)\times Sp(W)$-invariant.

In fact this definition generalizes the hyperdeterminant of boundary format
as introduced in Theorem 14.3.3 of \cite{GKZ}.

\begin{lemma0}
\rm If $A$ is degenerate then $D(A)=0$.
\end{lemma0}
\noindent {\bf Proof.}
There are $0\neq v\in V^*$ and $0\neq i\in I^*$ such that $A(v\otimes i)=0$.
 Hence $v\otimes S^{n+1}i\in V^*\otimes S^{n+1}I^*$  goes to zero
 under the dual of (\ref{magic}). \qed

\vspace{3mm}

If $A$ is nondegenerate we get $D(A)\neq 0$ only in the case $k=1$ and in
general
it can happen that $D(A)=0$, because the codimension of $N$ is $k$.
Our main technical result is the following

\begin{thm0}
\label{nonzero} If $A$ defines an instanton (that is $A$
belongs to ${\cal Q}^0$) then
$D(A)\neq 0$.
\end{thm0}

\noindent {\bf Proof.} From (\ref{definstanton})
we get the exact sequence
\refstepcounter{equation}
 \[ (\theequation) \label{ker} \hspace{8mm}
  0\rig{}K\rig{}W\otimes\OO\rig{}I\otimes\OO(1)\rig{}0. \]
The $(n+1)$-th wedge power twisted by $\OO(-n)$ gives the exact sequence
\[0\rig{}\wedge^{n+1}K(-n)\rig{}\wedge^{n+1}W(-n)\rig{}\ldots\]
\[\ldots\rig{}
\wedge^2W\otimes S^{n-1}I(-1)\rig{}W\otimes S^nI\rig{}S^{n+1}I(1)\rig{}0\]
where the $H^0$ of the last morphism corresponds to
$\Delta_A$ in (\ref{magic}).
Taking cohomology, it is enough to prove
\refstepcounter{equation}
 \[ (\theequation) \label{cla} \hspace{8mm}
H^n(\wedge^{n+1}K(-n))=0.  \]
The $(n+1)$-th wedge power twisted by $\OO(-n)$ of the sequence
\[0\rig{}I^*\otimes\OO(-1)\rig{}K\rig{}E\rig{}0\]
gives the sequence
\[ \begin{array}{c}0\rig{} S^{n+1}I^* \otimes K(-2n-1)\rig{}\ldots\rig{}
\wedge^{n-1}K\otimes S^{2}I^*(-n-2)\rig{}  \\
\rig{}\wedge^nK\otimes I^*(-n-1)  \rig{}\wedge^{n+1}K(-n)\rig{}\wedge^{n+1}E(-n
)
\rig{}0. \end{array}\]

In order to prove (\ref{cla}), taking cohomology, we need
$H^{n+i}(\wedge^{n-i}K(-n-i-1))=0$ for $i=0,\ldots, n$
and $H^n(\wedge^{n+1}E(-n))=0$.
The first group of vanishing is easily obtained by taking suitable
wedge powers of (\ref{ker}).
The crucial point to get the last vanishing is the isomorphism
 $\wedge^{n+1}E\simeq\wedge^{n-1}E$,
that it is true because $E$ is a rank $2n$ vector bundle with $c_1=0$. From
the sequence
\[0\rig{} S^{n-1}I^*(-2n-1)\rig{}S^{n-2}I^*\otimes K(-2n)
\rig{}\ldots\]
\[\ldots\rig{}\wedge^{n-1}K(-n)\rig{}
\wedge^{n-1}E(-n)\rig{}0,\]
in order to prove $H^n(\wedge^{n-1}E(-n))=0$, we only need
to see that \[H^{n+i}(\wedge^{n-1-i}K(-n-i))=0\quad\hbox{ for }
i=0,\ldots,n\]
 which
follows by using the exact sequence (\ref{ker}) exactly as above. \qed

\begin{corollary0}
\label{affine}
 $MI_{\PP^{2n+1}}(k)$ is affine.
\end{corollary0}

\noindent {\bf Proof.} By the \thmref{nonzero} we get that
${\cal Q}\setminus{N}={\cal Q}^0={\cal Q}\setminus \{D=0\}$ is affine.
It follows that
$MI_{\PP^{2n+1}}(k)$ is affine too, because it is the quotient of an affine
variety by a reductive group, see e.g. \cite{PV}, section 4.4. \qed

\begin{remark0}
\rm The invariant $D$ is meaningful even in the case $n=0$.
In this case it corresponds to the usual determinant of the map
$\CC^{2k}\to\CC^2\otimes\CC^k$.
For example for $n=0$ and $k=2$ the degenerate $2\times 2\times 4$ matrices
fill a variety of codimension $2$  and degree $12$ (\cite{BS})
in $\PP^{15}$ whose ideal
is generated by one quartic
(which is our invariant $D$), $10$ sextics and one octic.
We remark that the case $2\times 2\times 3$ is of boundary format.
The case $2\times 2\times 5$ is interesting. Here degenerate matrices fill a
variety of codimension $3$
and degree 20,
and its ideal is generated (at least) by $5$ quartics, $50$ sextics and $12$
octics.
The $5$ quartics
define a variety of codimension $2$ and degree $10$. Hence in this case no
analog of the invariant $D$ can exist.
\end{remark0}

\section{Instanton bundles with structure group $GL(2n)$} \indent

\begin{defn0}
\label{def:gral}
A $k$-instanton bundle $E$ on $\PP^{2n+1}$
is the cohomology bundle of a monad
\refstepcounter{equation}
 \[ (\theequation) \label{definstantong} \hspace{8mm}
 K\otimes\OO(-1)\rig{A} W\otimes\OO\rig{B} I\otimes\OO(1) \]
where $W$ is a  complex vector space of
dimension $2n+2k$ and $I,K$ are  complex vector spaces of
dimension $k$.
\end{defn0}

Notice that $E$ is not necessarily symplectic and that
this notion is a true generalization of the above one only for $n\ge 2$
because all rank $2$ bundles on $\PP^3$ with $c_1=0$ are symplectic.

 Let $(A,B) \in Hom(K \otimes V^*,W) \times Hom(W, I \otimes V)$ defining $E$.
The monad condition is now
equivalent
to the fact that
the matrices $A$ and $B$ are both nondegenerate and $B\cdot A=0$.

$Hom(K \otimes V^*,W) \times Hom(W, I \otimes V)$ contains the subvariety
${\cal P}$
 given by pairs of matrices $(A,B)$
for which the sequence (\ref{definstantong}) is a complex,
 that is, such that $B\cdot A=0$. $GL(I)\times GL(K)\times GL(W)$
acts on ${\cal P}$  by $(a,b,c)\cdot (A,B)=
(cAb, aBc^{-1})$.

As in the previous section we can see that there is a natural 1:1
correspondence between

\noindent {i)} isomorphism classes of instanton bundles, and

\noindent{ii)} orbits of $GL(I)\times GL(K)\times GL(W)$ on the open
subvariety ${{\cal P}^0}$ of ${\cal P}$ given by pairs of nondegenerate
matrices.

Moreover, if we denote by $H$ the quotient of
 $GL(I)\times GL(K)\times GL(W)$
 by $ (\lambda\cdot id,{\lambda}^{-1} \cdot id,\lambda \cdot id)$, then $H$
acts
 freely on ${{\cal P}^0}$.
 In particular, all points of ${{\cal P}^0}$ are stable for the action
 of $GL(I)\times GL(K) \times GL(W)$.

\begin{defn0}
\label{gitquoti}
The GIT-quotient ${{\cal P}^0}/{GL(I)\times GL(K) \times GL(W)}$
is denoted by \break\noindent $I_{\PP^{2n+1}}(k)$ and it is called
the moduli space of $k$-instanton bundles on $\PP^{2n+1}$.
It is a geometric quotient.
\end{defn0}

 $I_{\PP^{2n+1}}(k)$ coincides for $n=1,2$ with the open subset
 ${\cal I}_{\PP^{2n+1}}(k)$ of the Maruyama scheme of stable bundles
  on $\PP^{2n+1}$ of rank $2n$ and Chern polynomial $\frac{1}{(1-t^2)^k}$
   which are instanton bundles.
 For $n\ge 3$ we can  say that ${\cal I}_{\PP^{2n+1}}(k)$ is an open subset
 of $I_{\PP^{2n+1}}(k)$.
We remark that $MI_{\PP^{3}}(k)=I_{\PP^{3}}(k)$.
${\cal I}_{\PP^{2n+1}}(k)$ is known to be singular for
$n\ge 2$ and $k\ge 3$ (see \cite{MO}) and reducible for
 $n\ge 4$ (see \cite{AO00}).

\begin{defn0}
\label{def:bimagic}
Let $(A,B) \in Hom(K \otimes V^*,W) \times Hom(W, I \otimes V)$.
 We define
  \[ \tilde{D}(A,B):=\det S(A)\cdot\det R(B) \]
  where $\det$ denotes the usual determinant and $S(A)$, $R(B)$ are the
  following morphisms
  \[ \begin{array}{l} S(A): S^{n+1}K \otimes V^* \rightarrow S^{n}K \otimes W,
  \\ R(B): S^{n}I \otimes W \rightarrow S^{n+1}I \otimes V, \end{array}\]
  induced by $A$ and $B$ respectively, as in  \defref{def:magic}.
\end{defn0}

\begin{thm0}
\label{nonzero2} If $(A,B)$ defines an instanton (that is
$(A,B)$ belongs to ${\cal P}^{0}$) then
$\tilde{D}(A,B)\neq 0$.
\end{thm0}

\noindent {\bf Proof.} First of all we will see that
$\det S(A) \neq 0$. From (\ref{definstantong})
we get the exact sequence
\refstepcounter{equation}
 \[ (\theequation) \label{ker2} \hspace{8mm}
  0\rig{}K \otimes \OO(-1) \rig{}W\otimes\OO\rig{} Q \rig{}0. \]
The $(n+1)$-th wedge power twisted by $\OO(-n-2)$ gives the exact sequence
\[ \begin{array}{c}0\rig{}S^{n+1}K\otimes \OO(-2n-3)
\rig{}S^{n}K \otimes W \otimes \OO(-2n-2)\rig{}\ldots \\ \rig{}
\wedge^{n+1}W\otimes \OO(-n-2)\rig{}\wedge^{n+1}Q(-n-2) \rig{}0 \end{array} \]
where the $H^{2n+1}$ of the first morphism corresponds to $S(A)$. Hence,
taking cohomology, it is enough to prove
\[ H^n(\wedge^{n+1}Q(-n-2))=0.  \]
This is shown considering the $(n+1)$-wedge sequence of the exact sequence
\[0\rig{}E \rig{} Q\rig{}I\otimes\OO(1)\rig{}0\]
and arguing as in the proof of \thmref{nonzero}.

In order to prove $\det R(B) \neq 0$ we proceed exactly as in
 \thmref{nonzero}
and we leave the details to the reader.
 \qed

\begin{corollary0}
\label{affine2}
 $I_{\PP^{2n+1}}(k)$ is affine.
\end{corollary0}

\noindent {\bf Proof.}  First of all notice that given $(A,B) \in {\cal P}$,
if $A$ or $B$ is degenerate then $\det S(A)\cdot\det R(B)=0$. Hence, by
 \thmref{nonzero2} we get that
${\cal P}^0={\cal P}\setminus \{ \tilde{D}=0 \}$ is affine.
Therefore, by \cite{PV} section 4.4,
$I_{\PP^{2n+1}}(k)$ is affine too. \qed

\label{biblio}

{\bf Authors' addresses}
\vskip0.5truecm
\noindent
\hbox  {\hsize15truecm
\vbox  {
\hsize6truecm
\noindent Laura Costa\hfill \break
Dept. Algebra y Geometria \hfill \break
Universitat de Barcelona\hfill \break
Gran Via, 585\hfill \break
08007 Barcelona\hfill \break
SPAIN\hfill\break
costa@mat.ub.es}
\hfill
\vbox {
\hsize7truecm
\noindent Giorgio Ottaviani\hfill \break
Dipartimento di Matematica
"U. Dini" \hfill \break
Universit\`a di Firenze\hfill \break
viale Morgagni 67/A\hfill \break
I 50134 FIRENZE\hfill \break
ITALY\hfill\break
ottavian@math.unifi.it
}
}

\end{document}